
\def\LaTeX{%
  \let\Begin\begin
  \let\End\end
  \let\salta\relax
  \let\finqui\relax
  \let\futuro\relax}

\def\UK{\def\our{our}\let\sz s}
\def\USA{\def\our{or}\let\sz z}

\UK 



\LaTeX

\USA


\salta

\documentclass[twoside,12pt]{article}
\setlength{\textheight}{24cm}
\setlength{\textwidth}{16cm}
\setlength{\oddsidemargin}{2mm}
\setlength{\evensidemargin}{2mm}
\setlength{\topmargin}{-15mm}
\parskip2mm


\usepackage[raggedright]{titlesec}
\titleformat{\section}{\normalfont\large\bfseries}{\thesection}{1em}{}
\titleformat{\subsection}{\normalfont\bfseries}{\thesubsection}{1em}{}


\usepackage[english]{babel}
\usepackage[utf8]{inputenc}
\usepackage[usenames,dvipsnames]{color}
\usepackage{amsmath}
\usepackage{amsthm}
\usepackage{amssymb,bbm}
\usepackage[mathcal]{euscript}

\usepackage{cite}
\usepackage{hyperref}
\usepackage{enumitem}

\usepackage[ulem=normalem,draft]{changes}
%
%

%
 
\definecolor{viola}{rgb}{0.3,0,0.7}
\definecolor{ciclamino}{rgb}{0.65,0,0.65}
\definecolor{blu}{rgb}{0,0,0.7}
\definecolor{rosso}{rgb}{0.85,0,0}
\definecolor{darkgreen}{rgb}{0,0.5,0}


\usepackage{hyperref} 
\makeatletter
\hypersetup{%
	colorlinks	=true,
	linkcolor	=blue,%
	citecolor	=darkgreen,%
	urlcolor	=blue
}




\bibliographystyle{abbrv}


\numberwithin{equation}{section}
\newtheorem{theorem}{Theorem}[section]

\newtheorem{proposition}[theorem]{Proposition}
\newtheorem{lemma}[theorem]{Lemma}

\finqui

\def\Bcenter{\Begin{center}}
\def\Ecenter{\End{center}}
\let\non\nonumber




\def\step #1 \par{\medskip\noindent{\bf #1.}\quad}
\def\jstep #1: \par {\vspace{2mm}\noindent\underline{\sc #1 :}\par\nobreak\vspace{1mm}\noindent}




\def\multibold #1{\def\arg{#1}%
  \ifx\arg\pto \let\next\relax
  \else
  \def\next{\expandafter
    \def\csname #1#1\endcsname{{\boldsymbol #1}}%
    \multibold}%
  \fi \next}

\def\pto{.}

\def\multical #1{\def\arg{#1}%
  \ifx\arg\pto \let\next\relax
  \else
  \def\next{\expandafter
    \def\csname cal#1\endcsname{{\cal #1}}%
    \multical}%
  \fi \next}

\def\multigrass #1{\def\arg{#1}%
  \ifx\arg\pto \let\next\relax
  \else
  \def\next{\expandafter
    \def\csname grass#1\endcsname{{\mathbb #1}}%
    \multigrass}%
  \fi \next}

\def\multigrass #1{\def\arg{#1}%
  \ifx\arg\pto \let\next\relax
  \else
  \def\next{\expandafter
    \def\csname #1#1#1\endcsname{{\mathbb #1}}%
    \multigrass}%
  \fi \next}


\def\multimathop #1 {\def\arg{#1}%
  \ifx\arg\pto \let\next\relax
  \else
  \def\next{\expandafter
    \def\csname #1\endcsname{\mathop{\rm #1}\nolimits}%
    \multimathop}%
  \fi \next}

\multibold
qweryuiopasdfghjklzxcvbnmQWERTYUIOPASDFGHJKLZXCVBNM.  

\multical
QWERTYUIOPASDFGHJKLZXCVBNM.

\multigrass
QWERTYUIOPASDFGHJKLZXCVBNM.

\multimathop
diag dist div dom mean meas sign supp .


\def\Accorpa #1#2 #3 {\gdef #1{\eqref{#2}--\eqref{#3}}%
  \wlog{}\wlog{\string #1 -> #2 - #3}\wlog{}}


\def\esssup{\mathop{\rm ess\,sup}}

\def\<#1>{\mathopen\langle #1\mathclose\rangle}
\def\norma #1{\mathopen \| #1\mathclose \|}

\def\aeO{\checkmmode{a.e.\ in~$\Omega$}}
\def\aeQ{\checkmmode{a.e.\ in~$Q$}}

\def\Pi{\widehat\pi}

\def\phieps{\phi_\delta}
\def\mueps{\mu_\delta}

\def\iO{\int_\Omega}

\def\dt{\partial_t}

\def\checkmmode #1{\relax\ifmmode\hbox{#1}\else{#1}\fi}






\def\genspazio #1#2#3#4#5{#1^{#2}(#5,#4;#3)}
\def\spazio #1#2#3{\genspazio {#1}{#2}{#3}T0}

\def\L {\spazio L}
\def\H {\spazio H}

\def\C #1#2{C^{#1}([0,T];#2)}


\def\Lx #1{L^{#1}(\Omega)}
\def\Hx #1{H^{#1}(\Omega)}



\let\eps\varepsilon

\let\phi\varphi
\let\ph\phi

\let\TeXchi\chi                         
\newbox\chibox
\setbox0 \hbox{\mathsurround0pt $\TeXchi$}
\setbox\chibox \hbox{\raise\dp0 \box 0 }
\def\chi{\copy\chibox}



\def\0{{\boldsymbol {0} }}
\def\bet{{\eta}}

\def\Vp{{{H^1(\Omega)}^*}}


\def\CP0{(${\mathcal{CP}}_0$)}

\def\mobe{b_\delta}

\def\pheps{\phi_\delta}
\def\mueps{\mu_\delta}

\def\mob{b}
\def\mobe{\mob_\delta}

\DeclareMathAlphabet{\mathscr}{U}{mathc}{m}{it}

\Begin{document}


%
\title{
The anisotropic Cahn--Hilliard equation\\ with degenerate mobility:\\[0.5ex] 
Existence of weak solutions
}
\author{}
\date{}
\maketitle
\Bcenter
\vskip-1cm
{\large\sc Harald Garcke$^{(1)}$}\\
{\normalsize e-mail: {\tt\href{mailto:harald.garcke@ur.de}{harald.garcke@ur.de}}}
\\[0.25cm]
{\large\sc Patrik Knopf$^{(1)}$}\\
{\normalsize e-mail: {\tt\href{mailto:patrik.knopf@ur.de}{patrik.knopf@ur.de}}}
\\[0.25cm]
{\large\sc Andrea Signori$^{(2)}$}\\
{\normalsize e-mail: 
{\tt\href{mailto:andrea.signori@polimi.it}{andrea.signori@polimi.it}}
\\[0.25cm]
$^{(1)}$
{\small Faculty for Mathematics, University of Regensburg,}\\
{\small  93053 Regensburg, Germany}\\[.3cm]
$^{(2)}$
{\small Dipartimento di Matematica, Politecnico di Milano}\\
{\small via E. Bonardi 9, I-20133 Milano, Italy}\\ 
{\small Alexander von Humboldt Research Fellow}\\[.3cm]
}

\Ecenter

\medskip

\Begin{abstract}
\noindent 
This paper presents an existence result for the anisotropic Cahn--Hilliard equation characterized by a potentially concentration-dependent degenerate mobility taking into account an anisotropic energy. The model allows for the degeneracy of the mobility at specific concentration values, demonstrating that the solution remains within physically relevant bounds. The introduction of anisotropy leads to highly nonlinear terms making energy and entropy estimates rather involved. As the mobility degenerates in the pure phases, the degenerate Cahn--Hilliard equation describes  surface diffusion and is an important model to model solid-state dewetting (SSD) of thin films.
We show existence of weak solutions for the anisotropic degenerate Cahn--Hilliard equation by using suitable energy and entropy type estimates.



\End{abstract}

\thispagestyle{empty}

\newpage

\salta
\pagestyle{myheadings}
\newcommand\testopari{\sc Garcke  -- Knopf -- Signori}
\newcommand\testodispari{\sc The anisotropic degenerate Cahn--Hilliard equation}
\markboth{\testopari}{\testodispari}
\finqui
%

\section{Introduction}
\label{INTRO}
\setcounter{equation}{0}
The Cahn--Hilliard equation, originally introduced to model phase separation in binary alloys, has found widespread applications in describing a variety of phase transition processes in materials science, including solidification and the growth of thin solid films. In addition, there have been applications involving the  Cahn--Hilliard equation in diverse fields such as imaging sciences, two-phase flows and tumor growth. The classical Cahn--Hilliard equation 
\begin{align*}
    \dt \phi = b\Delta \mu,
    \quad 
    \mu =- \eps\Delta \phi + \frac 1\eps \psi'(\phi)
\end{align*}
describes the time evolution of a conserved order parameter $\phi$, typically representing the concentration difference of the two components of a binary mixture.  The equation can be derived as the scaled $H^{-1}$-type gradient flow (see \cite{TaylorC94, BDGP}) of the Ginzburg--Landau energy
\begin{align*}
	{\cal E}_{GL} (\phi)
	= \frac \eps2\iO |\nabla \phi|^2 + \frac 1 \eps \iO \psi(\phi),
\end{align*}
where the parameter $\eps$ represents a small positive constant related to the width of the diffuse interface between different phases, $\psi$ is a suitable double-well 
free energy density and $b$ is a constant mobility parameter.
For the above Cahn--Hilliard equation, first well-posedness results can be found in \cite{ElliottZ86} and since then many results have been obtained. We refer to \cite{Elliott1989, Miranvillebook} for overviews. 

Already in the original derivation of the Cahn--Hilliard equation (see \cite{Cahn1961spinodal}) a concentration dependent mobility was considered.
This leads to the diffusion equation
$$\dt \phi - \div (\mob(\phi) \nabla \mu) = 0.$$
Taking a mobility $\mob(\phi)=(1-\phi^2)/\eps$,  which degenerates at the pure phases $\phi=\pm 1$,
Cahn, Elliott and Novick-Cohen \cite{Cahn1996cahn} used formally matched asymptotic expansions to show that in the sharp interface limit,
as $\eps$ tends to zero, the geometric evolution law motion by minus the Laplacian of mean curvature is obtained.
This law can be interpreted as motion by surface diffusion and can be viewed as a higher order analogue of motion by mean curvature.
In \cite{TaylorC94}, the relation between the degenerate Cahn--Hilliard model and motion by surface diffusion was discussed with the
help of an $ H^{-1}$-gradient flow perspective. In \cite{EG}, the existence of weak solutions to the Cahn--Hilliard equation was shown
for the case of the isotropic energy ${\cal E}_{GL}$ by means of suitable energy and entropy type estimates. We also refer
to \cite{gruen95} for an existence result for a related fourth order degenerate parabolic equation.

The energy ${\cal E}_{GL}$ is isotropic due to the fact that the energy density depends on $\nabla \phi $ only through its modulus and
not through its direction.
In many physical applications, especially those involving crystalline materials, the assumption of isotropic interfacial energy is often inadequate. 
In anisotropic systems, the energy density depends locally on the orientation of the interface.
An anisotropic interfacial energy becomes crucial as it significantly influences the dynamics by allowing the energy to depend on the orientation of the interface. This consideration is essential for accurately modeling processes like solidification or thin film growth, where directional properties are important. In particular, stationary states will differ as they are typically given as so-called Wulff shapes, see e.g., \cite{bgnreview,BDGP} and the references therein.

To describe such phenomena, we consider an {\it anisotropic} free energy ${\cal E}_{A, \psi}$ given by 
\begin{align*}
	{\cal E}_{A,\psi} (\phi)
	= \eps\iO A(\nabla \phi) + \frac 1 \eps \iO \psi(\phi) .
\end{align*}
The anisotropy density $\gamma: \RRR^d \to (0,+\infty)$, $\gamma \in C^2(\RRR^d\setminus\{0\}) \cap C^0(\RRR^d)$ is assumed to be homogeneous of degree one
and we set 
\begin{align*}
	A(p):= \frac 12 |\gamma(p)|^2
	\quad \text{for all $p \in \RRR^d$}.
\end{align*}
We point out that to any function $A$ that satisfies suitable assumptions (see~\ref{ass:anisotropy}), we can conversely find an associated anisotropy density $\gamma$, which is related to the corresponding sharp interface energy density, see \cite{TaylorC94, GKNZ} for details. In the following analysis, our focus will primarily be on the function $A$.

Taking an $H^{-1}$-gradient flow perspective, 
Taylor and Cahn  \cite{TaylorC94} also formulated anisotropic variants of
the Cahn--Hilliard equations and of surface diffusion. Specifically, they obtained
\begin{subequations}    
    \label{eq:first}
    \begin{alignat}{2}
    	\label{eq:1first}
    	& \dt \phi - \div (\mob(\phi) \nabla \mu) = 0 , 
    	\\
    	\label{eq:2first}
    	& \mu = - \eps\div(A'(\nabla \phi)) + \tfrac 1\eps \psi'(\phi). 
    \end{alignat}
\end{subequations}
Using the formally matched asymptotic expansion approach of Cahn, Elliott, and Novick-Cohen \cite{Cahn1996cahn}, it was shown in \cite{Ratz2006surface, Dziwnik17anisotropic, GKNZ} that, after an appropriate time rescaling, the above system converges to an anisotropic form of motion by surface diffusion.

In recent years, the mathematical community has made significant progress in addressing analytical challenges for isotropic systems. However, the development of a comprehensive analytical framework for the anisotropic Cahn--Hilliard equation is still in its early stages. Notably, we mention the recent works \cite{GKW, GKNZ}, and the references therein. Especially, in the case of a degenerate mobility no analytical results for the anisotropic Cahn--Hilliard equation were known so far.
In fact, with respect to the dependence of $A$ on $\nabla \phi$ it was stated in \cite{TaylorC94} that
``It is not clear whether making these assumptions of the dependence of various functions on $\nabla \phi$ creates insurmountable mathematical problems when $\nabla \phi$ is zero ...''. 
Since $A$ is anisotropic and positively homogeneous of degree two, the generic case is that 
$A'$ is not differentiable. In fact, $A'$ is differentiable only in the very simple case that $A$ is a positive definite quadratic form.
For this reason, the corresponding elliptic operator $\div(A'(\nabla \phi))$ is in general non-smooth. 
However, for  a weak solution theory it will be necessary to obtain $H^2$ control in space of solutions and it was not clear how to use entropy estimates to obtain $H^2$ control in the case of an anisotropic energy. By using an appropriate weak formulation and using a suitable regularity theory we will show in this paper
that an existence theory for
\eqref{eq:first} is possible. Let us mention an earlier result 
in \cite{dziwnik2019existence} for a degenerate mobility  and for a very particular anisotropic energy.
We also refer to the independent results in
\cite{elbar2025weaksolutionssharpinterface}, where the $H^2$ estimates developed in this paper (cf.~\eqref{EST:H2}) are used to study an anisotropic, inhomogeous Cahn--Hilliard equation with a disparate mobility.

Let us now present a detailed formulation of the problem. Let $\Omega \subset \RRR^d$, where $d \in \NNN$, be a $d$-dimensional flat torus, let $T > 0$ denote a fixed final time, and let $Q:= \Omega \times (0,T)$. The system we aim to analyze reads as
\begin{subequations}
    \label{Sys}
    \begin{alignat}{2}
    	\label{eq:1}
    	& \dt \phi - \div (\mob(\phi) \nabla \mu) = 0 \qquad&& \text{in $Q$},
    	\\
    	\label{eq:2}
    	& \mu = - \eps\div(A'(\nabla \phi)) + \tfrac 1\eps \psi'(\phi) \qquad&& \text{in $Q$},
    	\\
    	\label{eq:3}
    	& \phi\vert_{t=0} =\phi_0 
    	\qquad
    	 && \text{in $\Omega$}.
    \end{alignat}
\end{subequations}
In fact, the assumption that $\Omega$ is a $d$-dimensional flat torus means that 
$\Omega={(a_1,b_1)}\times ... \times {(a_d,b_d)}$ (with $a_i<b_i, \; i=1,...,d$),
where \eqref{eq:1} and \eqref{eq:2} are considered subject to periodic boundary conditions.
As it will not affect the upcoming analysis, we will set $\eps=1$ from this point onward.

The structure of this paper is as follows: 
After introducing the main assumptions and the notation used, we will state the main result which is the existence of a weak solution to system  \eqref{Sys} in Section 3. The proof of the existence of a weak solution is presented in Section 4, where an appropriate non-degenerate and non-singular approximation scheme is employed. With the help of suitable energy- and entropy type estimates, we are able to show a priori bounds for the approximate solutions. Here, in particular, it is crucial to control the $H^2$-norm in space, due to the highly nonlinear and non-smooth anisotropy function $A$. By means of the a priori estimates, we apply suitable compactness results, which make it possible to pass to the limit in the approximate equations. In an Appendix we eventually show an elementary result on a weak differentiability property of the composition of Lipschitz and Sobolev functions.

\section{Preliminaries}

\subsection{Assumptions}
Throughout this paper we will use the following assumptions.
\begin{enumerate}[label=$(\mathrm{A \arabic*})$, ref = $(\mathrm{A \arabic*})$,topsep=0ex]
\item \label{ass:Omega}
   The set $\Omega$ is a $d$-dimensional flat torus, with $d \in \NNN$,  and $T > 0$ denotes an arbitrary final time.   

\item \label{ass:anisotropy}
    For the anisotropy function $A$, we postulate  that $A\in C^{1,1}(\RRR^d)$     
    is positive on $\RRR^d \setminus\{0\}$ and positively two-homogeneous, that is,
    \begin{align*}
        A(\lambda p ) =\lambda^2 A(p), \quad \text{for all $\lambda >0$, $p \in \RRR^d$}.
    \end{align*}
    It then follows that there exist positive numbers $A_0, A_1, A_2$ with $0< A_0 < A_1$ such that, for all $p\in\RRR^d$, it holds
    \begin{align*}
        A_0 |p|^2
        \leq 
        A(p)
        \leq    
        A_1 |p|^2,
        \quad 
        |A'(p)|\leq A_2
        |p|.
    \end{align*}
    Besides, we require that the gradient of $A$, denoted by $A':\RRR^d\to \RRR^d$, is strongly monotone. This means that there exists a constant $c_A>0$ such that for all $p,q \in \RRR^d$, it holds 
    \begin{align*}
        (A'(p)-A'(q))\cdot (p-q)
        \geq c_A |p-q|^2.
    \end{align*}
    This readily entails that $A$ is strongly convex and, a fortiori, strictly convex.
   
\item \label{ass:pot}
    The potential $\psi$ can be decomposed as $\psi = \psi_1 + \psi_2$ with $\psi_2\in C^2([-1,1])$, and $\psi_1\in C^0([-1,1]) \cap C^2((-1,1))$ such that
    \begin{align}
        \label{ass:pot:growth}
        \psi_1''(r) = (1-r^2)^{-m}F(r) \quad \text{for all $r \in (-1,1)$},
    \end{align}
    for some $m\ge 1$, where 
    $F:[-1,1]\to [0,+\infty)$ is of class $C^1$. Without loss of generality, we assume that $\psi \ge 0$ on $[-1,1]$.

\item \label{ass:mob}
    The mobility function $\mob:\RRR \to [0,+\infty)$ is given by
    \begin{align}
        \label{ass:mob:formula}
        b(r)=
        \begin{cases}
            (1-r^2)^{m} B(r), \quad &\text{if $r\in [-1,1]$},\\
            0,\quad  &\text{if $r\in\RRR\setminus [-1,1]$},
        \end{cases}
    \end{align}
    where $m\ge 1$ is the constant from \ref{ass:pot}.
    Here, the function $B: [-1,1] \to (0,+\infty)$ is of class $C^1$ and there exist constants $0<B_*<B^*$ such that 
    \begin{align*}
            B_* \leq B(r) \leq B^*
            \quad \text{for all $r \in [-1,1]$}.
    \end{align*}
    Consequently, $b\vert_{[-1,1]}$ is continuously differentiable on $[-1,1]$ and positive on $(-1,1)$.
\end{enumerate}

As a consequence of \ref{ass:pot} and \ref{ass:mob}, we point out that
\begin{equation}
    \label{ID:BPSI}
    \mob \psi'' = FB  +  \mob\psi_2'' \quad\text{on $[-1,1]$},
\end{equation}
which entails that $\mob\psi''$ is bounded on $[-1,1]$.
We further define the function 
\begin{align}
    \label{def:Phi}
    \Phi:(-1,1) \to [0,+\infty),
    \quad
    r \mapsto \int_0^r \int_0^s \frac{1}{b(\tau)} \,\mathrm d\tau \,\mathrm ds .
\end{align}
Due to \ref{ass:mob}, $\Phi$ is nonnegative, twice continuously differentiable and satisfies
\begin{align}
    \label{prop:Phi}
    \Phi''(r) = \frac 1 {\mob(r)}
    \quad \text{for all $r \in (-1,1)$,}
    \qquad 
    \Phi(0)=\Phi'(0)=0.
\end{align}

\subsection{Notation} 
We further fix some notation that will be used in the remainder of this paper.

For any normed space $X$, we denote its norm by $\|\cdot\|_X$,
its dual space by $X^*$, and the duality pairing between $v \in X^*$ and $u\in X$ by $\langle v,u \rangle_X$.
Moreover, if $X$ is a Hilbert space, we write $(\cdot,\cdot)_X$ to denote its inner product.
    
In view of \ref{ass:Omega}, we write 
$L^p(\Omega)$ and $W^{k,p}(\Omega)$ (with $1 \leq p \leq \infty$ and $k \in \mathbb{N}_0$) to
denote the Lebesgue and Sobolev spaces of periodic functions on the flat torus $\Omega$.
Their corresponding norms are denoted by $\|\cdot\|_{L^p(\Omega)}$ and $\|\cdot\|_{W^{k,p}(\Omega)}$, respectively.
In the case $p = 2$, we use the notation $H^k(\Omega) = W^{k,2}(\Omega)$, and we identify $H^0(\Omega)$ with $L^2(\Omega)$.
    
For $k \in \mathbb{N}$, $1 \leq p \leq \infty$, any $T>0$ and any Banach space $X$, we write $L^p(0,T;X)$, $W^{k,p}(0,T;X)$ and $H^{k}(0,T;X) = W^{k,2}(0,T;X)$ to denote the Lebesgue and Sobolev spaces of functions with values in $X$. The standard norms on these spaces are denoted by $\|\cdot\|_{L^p(0,T;X)}$, $\|\cdot\|_{W^{k,p}(0,T;X)}$ and $\|\cdot\|_{H^k(0,T;X)}$, respectively. We further write $C^0([0,T];X)$ to denote the space of continuous functions mapping from $[0,T]$ to $X$. 

\pagebreak[4]

\section{Main result}

The main result of our paper is the existence of a weak solution to system \eqref{Sys}, which is stated in the following theorem. 

\begin{theorem}
\label{THM:MAIN}
Assume that \ref{ass:Omega}--\ref{ass:mob} are satisfied. Let $\varphi_0:\Omega\to\RRR$ be a prescribed initial datum, which satisfies
\begin{subequations}
    \label{ass:ini}
    \begin{align}
        \label{ass:init:data:1}
         & \ph_0 \in {\Hx1}, 
         \quad
         |\ph_0 | \leq 1 \quad \aeO,
         \\ &
         \label{ass:init:data:2}
         \psi(\ph_0) \in \Lx1,
         \quad  
         \Phi(\ph_0) \in \Lx1.
    \end{align}
\end{subequations}
Then, there exists a weak solution $(\ph,J,w)$ to the system \eqref{Sys}, which has the following properties:
\begin{enumerate}[label=\textnormal{(\roman*)}, leftmargin=*,topsep=0ex]
    \item \label{THM:MAIN:1} It holds
    \begin{align*}
    	&\ph  \in \H1 \Vp \cap \C0 {\Lx2}
        \cap \L\infty {\Hx1} \cap \L2 {{\Hx2}} ,
    	\\[0.5ex] 
    	&J \in \L2 {L^2(\Omega;\RRR^d)},
        \quad
        w \in \L2 {\Lx2},
        \quad 
        \psi'(\ph)  \in \L2 {\Lx2},
        \\[0.5ex]
        &A'(\nabla \ph) \in \L2 {\Hx1} ,
        \quad
        \mob(\ph) \in \L2 {\Hx1}
    \end{align*}
    as well as $|\varphi|\le 1$ a.e.~in $Q$.

    \item \label{THM:MAIN:2} The weak formulation
    \begin{subequations}
        \label{WF}
        \begin{align}
            \label{WF:1}
            \<\dt \ph, \zeta>_{\Hx1}
            &= 
            \iO J \cdot \nabla \zeta,
            \\
            \label{WF:2}
            \iO J \cdot \bet 
            &= \iO w \div \big( \mob(\ph)\bet \big) 
                - \mob(\ph) \psi''(\ph) 
                \nabla \ph \cdot \bet,
            \\
            \label{WF:3}
            \iO w \xi  
            &= \iO A'(\nabla \ph ) \cdot \nabla \xi,
        \end{align}
    \end{subequations}
    hold almost everywhere on $[0,T]$ for all test functions $\zeta,\xi\in H^1(\Omega)$ and $\bet\in H^1(\Omega;\RRR^d) \cap L^\infty(\Omega;\RRR^d)$. 
    Moreover, the initial condition $\ph\vert_{t=0}=\ph_0$ is fulfilled almost everywhere in $\Omega$.
    \end{enumerate}
\end{theorem}

The proof of this theorem will be presented in Section~\ref{SECT:PROOF}.

Theorem~\ref{THM:MAIN} extends the existence results established in \cite{GKNZ,GKW} for the anisotropic Cahn--Hilliard equation with a uniformly positive mobility function to degenerate mobilities. In contrast to the existence result developed in \cite{EG} for the isotropic Cahn--Hilliard equation (i.e., $A=\tfrac 12 |\cdot|^2$) with degenerate mobility in a bounded domain, our proof of Theorem~\ref{THM:MAIN} crucially relies on the framework of periodic boundary conditions.
As in the isotropic case, it remains an open problem whether higher regularity properties or even uniqueness of the weak solution can be obtained.

\pagebreak[4]

\section{Existence of a weak solution} \label{SECT:PROOF}

This section is devoted to the proof of our main result Theorem~\ref{THM:MAIN}.
Our strategy for the construction of a weak solution is to approximate the degenerate mobility function $b$, by a sequence $(b_\delta)_{\delta>0}$ of non-degenerate, uniformly positive mobility functions. Then, for any sufficiently small $\delta>0$, the theory developed in \cite{GKNZ} and \cite{GKW} can be applied and provides the existence of a weak solution as well as regularity properties for the approximate problem. Adapting the approach in \cite{EG}, we derive suitable uniform estimates (with respect to $\delta$) for the approximate solutions. Via compactness arguments, we eventually pass to the limit $\delta\to 0$, and show that the approximate solutions converge to a weak solution of the original system \eqref{Sys} as stated in the theorem. 

\subsection{Approximation scheme}
We start by introducing our approximate problem. 
Let $\delta\in (0,1)$ be an arbitrary, prescribed real number.
We define the approximated mobility functions as
\begin{align}
	\label{mob:approx}
    \mobe:\RRR \to \RRR, \quad 
	\mobe(r) =
	\begin{cases}
	\mob(-1+\delta), \quad & \text{if $r \leq \delta-1$,}
	\\
	\mob(r), \quad & \text{if $|r| < 1-\delta$,}
	\\
	\mob(1-\delta), \quad & \text{if $r \geq 1-\delta$.}
	\end{cases}
\end{align}
By this construction, due to \ref{ass:anisotropy}, $\mobe$ is differentiable on $\RRR\setminus\{-1+\delta,1-\delta\}$, and its derivative is given by
\begin{align*}
    \mobe'(r) =
    \begin{cases}
        0, &\text{if $|r|>1-\delta$},\\
        b'(r), &\text{if $|r|<1-\delta$},
    \end{cases}
\end{align*}
and it holds
\begin{align}
    \label{EST:mobLip}
    \big| \mobe'(r) \big| \le \norma{b'}_{L^\infty([-1,1])} \le C
    \quad\text{for all $r\in\RRR$}.
\end{align}
This implies that $\mobe$ is Lipschitz continuous with a Lipschitz constant independent of $\delta$. 

We next introduce an approximation of the function $\Phi$ defined in \eqref{def:Phi} as 
\begin{align}
    \label{def:PhiE}
    \Phi_\delta:\RRR \to [0,+\infty),
    \quad
    r \mapsto \int_0^r \int_0^s \frac{1}{\mobe(\tau)} \,\mathrm d\tau \,\mathrm ds .
\end{align}
Due to the properties of $\mobe$, it is straightforward to check that $\Phi_\delta$ is twice continuously differentiable and satisfies
\begin{align}
    \label{COND:PHIE:1}
    \left\{
    \begin{aligned}
    \Phi_\delta''(r) &= \frac 1 {\mobe(r)},
    &&\quad \text{for all $r \in \RRR$,}
    \\
    \Phi_\delta(r) &= \Phi(r) \phantom{\Big|},
    &&\quad \text{for all $r \in [-1+\delta,1-\delta]$,}
    \\[1ex]
    \Phi_\delta(0) &= \Phi_\delta'(0)=0.
    \end{aligned}
    \right.
\end{align}
Moreover, there exists a constant $c_\delta \ge 0$ such that
\begin{equation}
    \label{COND:PHIE:2}
    |\Phi_\delta''(r)| \le c_\delta
    \quad\text{and}\quad
    |\Phi_\delta'(r)| \le c_\delta |r|
    \quad\text{for all $r\in\RRR$.}
\end{equation}

We further approximate the singular potential $\psi$ by regular potentials $\psi_{\delta} := \psi_{1,\delta} + \psi_2$, where
$\psi_{1,\delta}:\RRR \to\RRR$ is given by
\begin{equation}
    \label{DEF:POT:APPROX}
    \psi_{1,\delta}(r) =
    \begin{cases}
	\displaystyle \sum_{k=0}^2 \frac{\psi_1^{(k)}(-1+\delta)}{k!}[r - (-1+\delta)]^k, \quad & \text{if $r \leq -1+\delta$,}
	\\
	\psi_1(r) \phantom{\Bigg|}, \quad & \text{if $|r| < 1-\delta$,}
	\\
	\displaystyle \sum_{k=0}^2 \frac{\psi_1^{(k)}(1-\delta)}{k!}[r - (1-\delta)]^k, \quad & \text{if $r \geq 1-\delta$.}
	\end{cases}
\end{equation}
In this way, $\psi_{1,\delta}\in C^2(\RRR)$ with
\begin{align}
	\label{pot:approx}
	\psi''_{1,\delta}(r) :=
	\begin{cases}
	\psi_1''(-1+\delta), \quad & \text{if $r \leq -1+\delta$,}
	\\
	\psi_1''(r), \quad & \text{if $|r| < 1-\delta$,}
	\\
	\psi_1''(1-\delta), \quad & \text{if $r \geq 1-\delta$,}
	\end{cases}
\end{align}
and it holds $\psi_{1,\delta}(0)=\psi_{1}(0)$, and $\psi_{1,\delta}'(0)=\psi_{1}'(0)$.

Lastly, we also approximate our initial datum by defining
\begin{equation*}
    \phi_{0,\delta} := (1-\delta) \phi_0,
\end{equation*}
which ensures that $|\phi_{0,\delta}| \le 1-\delta$ almost everywhere in $\Omega$.

We now introduce our approximate problem, which reads as
\begin{subequations}
    \label{Syseps}
    \begin{alignat}{2}
    	\label{eq:eps:1}
    	& \dt \pheps - \div (\mobe(\pheps) \nabla \mueps) = 0 \qquad&& \text{in $Q$},
    	\\
    	\label{eq:eps:2}
    	& \mueps = -  \div(A'(\nabla \pheps)) +  \psi_\delta'(\pheps) \qquad&& \text{in $Q$},
    	\\
    	\label{eq:eps:3}
    	& \pheps|_{t=0} = \phi_{0,\delta}
    	\qquad
    	 && \text{in $\Omega$}.
    \end{alignat}
\end{subequations}
The well-posedness of system \eqref{Syseps} as well as further regularity properties are ensured by the following proposition.

\begin{proposition}
\label{PROP:REGMOB:GKW}
Suppose that \ref{ass:Omega} and \ref{ass:anisotropy} are fulfilled and that the initial datum fulfills \eqref{ass:ini}.
We further assume that the mobility function $\mob$ is continuous and positive, that is, it fulfills
\begin{align}
    \label{mob:reg}
    0<m_*\leq \mob(r) \leq M^* \quad \text{for all $r \in \RRR,$}
\end{align}
for suitable positive constants  $m_*$ and $M^*$.
Moreover, we assume that the potential $\psi$ is continuously differentiable
and there exist nonnegative constants $c_1$, $c_2$ and $c_3$ as well as an exponent $\sigma \in [2,6)$ such that 
\begin{align}
    \label{pot:reg}
    - c_1 \leq \psi (r )
    \leq c_2 (|r|^\sigma+1),
    \quad 
    |\psi' (r )|
    \leq c_3 (|r|^{\sigma-1}+1)
    \quad \text{for all $r \in \RRR.$}
\end{align}
Then, there exists a weak solution $(\ph,\mu)$ to system \eqref{Sys} with the following properties:
\begin{enumerate}[label=\textnormal{(\roman*)}, leftmargin=*, topsep=0ex]
    \item The functions $\ph$ and $\mu$ have the regularities
    \begin{align}
    \label{REG:DNG}
    \left\{
    \begin{aligned}
    	&\ph  \in \H1 \Vp \cap \C0 {\Lx2}
        \\
        &\qquad \cap \L\infty {\Hx1} \cap \L2 {{\Hx2}} ,
    	\\[0.5ex] 
    	&\mu  \in \L2 {\Hx1},
        \quad 
        \psi'(\ph)  \in \L2 {\Lx2},
        \\[0.5ex]
        &A'(\nabla \ph(t)) \in {\Hx1} \quad \text{for almost all $t \in [0,T].$}
    \end{aligned}
    \right.
\end{align}

    \item It holds
    \begin{subequations}
    \label{WF:NDG}
        \begin{align}
            \label{WF:NDG:1}
        	\<\dt \ph, \zeta>_{\Hx1}
        	&= 
        	- \iO \mob(\ph) \nabla \mu \cdot \nabla \zeta
         \end{align}
    almost everywhere on $[0,T]$ for all $ \zeta \in {\Hx1}$, and
        \begin{align}
            \label{WF:NDG:2}
            \mu = -\div\big(A'(\nabla\phi)\big)  + \psi'(\phi)
            \quad\aeQ.
        \end{align}
    \end{subequations}    
    Moreover, the initial condition $\ph\vert_{t=0}=\ph_0$ is fulfilled almost everywhere in $\Omega$.

    \item The weak energy dissipation law
    \begin{equation}
        \label{IEQ:EN}
        {\cal E}_{A,\psi} \big(\phi(t)\big)
        + \frac 12 \int_0^t \iO b(\phi) |\nabla\mu|^2
        \le
        {\cal E}_{A,\psi} \big(\phi_0\big)
    \end{equation}
    holds for almost all $t\in [0,T]$.
    
\end{enumerate}
\end{proposition}

For a proof of this lemma, we refer to \cite[Theorem 4.3]{GKNZ} and \cite[Proposition 2.2 and Theorem 2.8]{GKW}. To be precise, the aforementioned results in \cite{GKNZ} and \cite{GKW} were established in the case of a bounded domain $\Omega$, where Neumann boundary conditions were imposed on $\ph$ and $\mu$. However, it is easy to transfer these results to the situation, where $\Omega$ is a flat torus. In fact, some parts of the mathematical analysis (e.g., the $H^2$-regularity theory) become a lot easier in the periodic setting since boundary effects can be ignored.

In view of \eqref{mob:approx}, it is clear that $\mobe$ fulfills assumption \eqref{mob:reg}. Moreover, it follows from \eqref{DEF:POT:APPROX}, \eqref{pot:approx} and the assumptions on $\psi_2$ that the conditions in \eqref{pot:reg} are fulfilled by $\psi_\delta$ with $\sigma=2$.
Hence, applying Proposition~\ref{PROP:REGMOB:GKW} to our approximate problem \eqref{Syseps}, we directly infer
the existence of a weak solution $(\pheps,\mueps)$ with 
\begin{align}
    \label{REG:DNG:E}
    \left\{
    \begin{aligned}
    	&\pheps  \in \H1 \Vp \cap \C0 {\Lx2}
        \\
        &\qquad \cap \L\infty {\Hx1} \cap \L2 {{\Hx2}} ,
    	\\[0.5ex] 
    	&\mu_\delta  \in \L2 {\Hx1},
        \quad 
        \psi_\delta'(\pheps)  \in \L2 {\Lx2},
        \\[0.5ex]
        &A'(\nabla \pheps(t)) \in {\Hx1} \quad \text{for almost all $t \in [0,T].$}
    \end{aligned}
    \right.
\end{align}

To prove that the approximate solutions $(\pheps,\mueps)$ converge to a weak solution of the original problem \eqref{Sys}, we need to show that they fulfill a variational formulation that is compatible with \eqref{WF}. To this end, we introduce the auxiliary functions
\begin{align}
    \label{DEF:JE}
    J_\delta &:= - \mobe(\pheps) \nabla \mueps,
    \\
    \label{DEF:WE}
    w_\delta &:= - \div \big(A'\big(\nabla\pheps\big)\big),
\end{align}
and claim that the weak formulations
\begin{subequations}
    \label{WF:NDG:ALT}
    \begin{align}
        \label{WF:NDG:ALT:1}
        \<\dt \pheps, \zeta>_{\Hx1}
        &= 
        \iO J_\delta \cdot \nabla \zeta,
        \\
        \label{WF:NDG:ALT:2}
        \iO J_\delta \cdot \bet 
        &= \iO w_\delta \div \big( \mobe(\pheps)\bet \big) 
            - \mobe(\pheps) \psi_\delta''(\pheps) 
            \nabla \pheps \cdot \bet,
        \\
        \label{WF:NDG:ALT:3}
        \iO w_\delta \xi  
        &= \iO A'(\nabla \pheps ) \cdot \nabla \xi,
    \end{align}
\end{subequations}
are fulfilled almost everywhere on $[0,T]$ for all test functions $\zeta,\xi\in H^1(\Omega)$ and $\bet\in H^1(\Omega;\RRR^d) \cap L^\infty(\Omega;\RRR^d)$. 

In fact, \eqref{WF:NDG:ALT:1} directly follows from the weak formulation \eqref{WF:NDG:1} written for $\mobe$, $\psi_\delta$ and $(\pheps,\mu_\delta)$ along with the definition of $J_\delta$ in \eqref{DEF:JE}. The identity \eqref{WF:NDG:ALT:1} is simply the weak formulation of \eqref{DEF:WE}.

To derive \eqref{WF:NDG:ALT:2}, we first recall \eqref{EST:mobLip}, which entails that $\mobe$ is
Lipschitz continuous with a Lipschitz constant independent of $\delta$. Moreover, it holds
\begin{align}
    \label{EST:MOBE}
    \norma{\mobe(\pheps) }_{L^\infty(Q)}
    \le \norma{b}_{L^\infty(\RRR)}
\end{align}
for almost all $t\in [0,T]$.
Thus, Lemma~\ref{LEM:LIPSCHITZ} implies that, for almost all $t\in [0,T]$, it holds
\begin{align}
    \label{EST:DMOBE}
    \mobe\big(\pheps(t)\big) \in H^1(\Omega)
    \quad\text{with}\quad   
    \big\| \mobe\big(\pheps(t)\big) \big\|_{H^1(\Omega)}
    &\le C \norma{\nabla\pheps(t)}_{H^1(\Omega)}.
\end{align}
Note that according to Lemma~\ref{LEM:LIPSCHITZ}, the constant $C$ on the right-hand side can be chosen independent of $\delta$ as the Lipschitz constant of $\mobe$ is independent of $\delta$.
Moreover, using the chain rule formula for the combination of a scalar Lipschitz function with a Sobolev function (see, e.g., \cite[Corollary~3.2]{Ziemer}), we infer
\begin{align}
    \label{CR:MOBE}
    \nabla\big(\mobe(\pheps)\big)(x,t) = 
    \mobe'\big(\pheps(x,t)\big) \nabla\pheps(x,t) 
    \quad
    \text{for almost all $(x,t)\in Q$}. 
\end{align}
Let now $\bet\in H^1(\Omega;\RRR^d) \cap L^\infty(\Omega;\RRR^d)$ be arbitrary. Then, the function
\begin{equation*}
    \vartheta := 
    \div\big( \mobe(\pheps) \bet \big) 
    = \mobe'(\pheps) \nabla\pheps \cdot \eta
        + \mobe(\pheps) \div(\bet)
\end{equation*}
belongs to $L^2(\Omega)$ almost everywhere on $[0,T]$.
Hence, using \eqref{DEF:JE} as well as \eqref{WF:NDG:2} written for $\psi_\delta$ and $(\pheps,\mu_\delta)$, we infer via integration by parts that
\begin{align}
    \non 
    &\iO J_\delta \cdot \bet 
    = - \iO \nabla \mu_\delta \cdot \big( \mobe(\pheps) \bet \big)
    = \iO \mu_\delta \vartheta
    \\ \non
    &\quad = \iO \big[ -\div\big( A'(\nabla\phi_\delta) \big) + \psi_\delta'\big(\pheps\big) \big]
        \vartheta
    \\ \non
    &\quad
    = \iO \big[ w_\delta + \psi_\delta'\big(\pheps\big)\big] 
        \div\big( \mobe(\pheps) \bet \big) 
    \\
    &\quad
    = \iO w_\delta \div \big( \mobe(\pheps)\bet \big) 
            - \mobe(\pheps) \psi_\delta''(\pheps) 
            \nabla \pheps \cdot \bet
            \label{DEF:JE:wf}
\end{align}
almost everywhere on $[0,T]$. This proves that \eqref{WF:NDG:ALT:2} holds for all $\bet\in H^1(\Omega;\RRR^d) \cap L^\infty(\Omega;\RRR^d)$.

\subsection{Uniform bounds}

In the following computations, the letter $C$ will denote general nonnegative constants that may depend on the initial data and the quantities introduced in \ref{ass:Omega}--\ref{ass:mob}, but is independent of $\delta$.

Since $\psi$ is continuous on $[-1,1]$ (cf.~\ref{ass:pot}) and $\phi_{0,\delta} \to \phi_0$ pointwise almost everywhere in $\Omega$ as $\delta\to 0$, it also holds that
\begin{equation*}
    \psi(\phi_{0,\delta}) \to \psi(\phi_{0})
    \quad\text{\aeO\ as $\delta\to 0$.}
\end{equation*}
As we further have $|\psi(\phi_{0,\delta})| \le \norma{\psi}_{L^\infty([-1,1])}$ almost everywhere in $\Omega$, Lebesgue's dominated convergence theorem yields
\begin{equation}
    \iO \psi(\phi_{0,\delta}) \to \iO \psi(\phi_{0})
    \quad\text{as $\delta\to 0$.}
\end{equation}
Since $\psi$ is nonnegative (cf.~\ref{ass:pot}) and $A$ is positively two-homogeneous (cf.~\ref{ass:anisotropy}), we infer that
\begin{equation*}
    {\cal E}_{A,\psi} (\phi_{0,\delta})
    = (1-\delta)^2 \iO A(\nabla \phi_0) + \iO \psi(\phi_{0,\delta})
    \le {\cal E}_{A,\psi} (\phi_{0}) + 1,
\end{equation*}
where here and in the following, we assume $\delta$ to be chosen sufficiently small.
Recalling further that $|\phi_{0,\delta}| \le 1-\delta$ almost everywhere in $\Omega$,
we conclude by means of the energy inequality \eqref{IEQ:EN} written for $\psi_\delta$ that
\begin{align}
	\label{E:eps}
	{\cal E}_{A,\psi_\delta}(\ph (t))
	\leq {\cal E}_{A,\psi_\delta} (\phi_{0,\delta})
    = {\cal E}_{A,\psi} (\phi_{0,\delta})
    \leq {\cal E}_{A,\psi} (\phi_0) + 1
\end{align}
for almost all $t\in [0,T]$. Invoking assumption~\ref{ass:anisotropy}, we directly infer
\begin{align}
    \label{UNI:1}
    \norma{\nabla\pheps}_{\L\infty {{L^2(\Omega;\RRR^d)}}}^2
        + \int_0^T \iO \mobe(\pheps) |\nabla\mueps|^2
    \le C.
\end{align}
Testing \eqref{WF:NDG:ALT:1} with $\zeta\equiv 1$, we infer the conservation of mass, that is
\begin{equation}
    \label{UNI:2}
    \iO \phi_\delta (t) = \iO \phi_{0, \delta}
    \quad\text{for all $t\in[0,T]$.}
\end{equation}
Combining \eqref{UNI:1} and \eqref{UNI:2}, we employ Poincar\'e's inequality to obtain
\begin{align}
	\label{UNI:PHE}
	\norma{\pheps}_{\L\infty {\Hx1}}
	\leq C.
\end{align}
Moreover, recalling \ref{ass:mob} and \eqref{DEF:JE:wf}, we use \eqref{UNI:1} to further deduce
\begin{align}
    \label{UNI:JE}
    \norma{J_\delta}_{\L2 {{L^2(\Omega;\RRR^d)}}}^2 
    \le B_* \int_0^T \iO \mobe(\pheps) |\nabla\mueps|^2
    \le C.
\end{align}
Eventually, we conclude from \eqref{WF:NDG:ALT:1} that
\begin{align}
    \label{UNI:DTPHE}
    \norma{\partial_t \pheps}_{\L2 {H^1(\Omega)^*}} \le C.
\end{align}

Let us now proceed with a higher-order estimate inspired by \cite{EG}.  
We first notice that due to \eqref{COND:PHIE:2}, it holds $\Phi_\delta(\pheps) \in H^1(\Omega)$ almost everywhere on $[0,T]$. Testing \eqref{eq:eps:1} with $\Phi_\delta'(\pheps)$, and using \eqref{WF:NDG:2}, written for $(\pheps,\mueps)$ and $\psi_\delta$, as well as \eqref{WF:NDG:ALT:1}, we infer 
\begin{align*}
	&\frac{\mathrm d}{\mathrm dt} \iO \Phi_\delta(\pheps)
	= \langle \dt \pheps ,\Phi_\delta' (\pheps)\rangle_{H^1(\Omega)} 
	= - \iO \mobe(\pheps) \Phi_\delta''(\pheps) \nabla \mu_\delta \cdot \nabla \pheps
	\notag\\ 
    &\quad  = 
	\iO \mu_\delta \Delta \pheps
	= - \iO \div (A'(\nabla \pheps)) \Delta \pheps
	- \iO \psi_\delta''(\pheps)|\nabla \pheps|^2
\end{align*}
almost everywhere on $[0,T]$ (cf.~\cite{EG}).
Thus, upon integrating over time and rearranging some terms, 
we deduce that, for all $t\in [0,T]$,
\begin{align}
    \label{EQ:PHI}
	& \iO \Phi_\delta(\pheps(t))
	+ \int_0^t\iO \div (A'(\nabla \pheps)) \Delta \pheps
	+ \int_0^t\iO \psi_{1,\delta}''(\pheps)|\nabla \pheps|^2
	\notag\\ 
    &\quad  =
	\iO \Phi_\delta(\ph_{0,\delta})
	- \int_0^t\iO \psi_2''(\pheps)|\nabla \pheps|^2 .
\end{align}

We now want to estimate the second integral on the left-hand side of \eqref{EQ:PHI} from below.
In order to exploit the monotonicity of $A'$, we approximate some of the involved derivatives by difference quotients.
To this end, for any given scalar or vector-valued function $f$ defined on $\Omega$, we write $\partial_i^{+h}f$ and $\partial_i^{-h}f$ to denote the difference quotients at a given point $x\in \Omega$ by
\begin{align*}
    \partial_i^{+ h} f (x):=   \frac {f(x+h e_i) - f(x)}h,
    \quad\text{and}\quad
    \partial_i^{- h} f (x):=   \frac {f(x) - f(x-he_i)}h,
\end{align*}
respectively.
Recall that, due to \eqref{REG:DNG:E}, it holds that
\begin{align*}
    \pheps(t) \in H^2(\Omega)
    \quad\text{and}\quad
    A'\big(\pheps(t)\big) \in H^1(\Omega)
     \quad\text{for almost all $t\in [0,T]$}.
\end{align*}
Hence, we have
\begin{align}
    \label{CONV:PHEPSJI}
    \partial_j^{+h} \partial_i \pheps (t) \to \partial_j \partial_i \pheps (t) 
    \quad\text{in $L^2(\Omega)$ \,for all $i,j\in\{1,...,d\}$ as $h\to 0$},
\end{align}
as well as
\begin{equation}
    \label{CONV:PHEPSJJ}
    \sum_{j=1}^d \partial_j^{-h} \partial_j^{+h} \pheps(t) \to \Delta \pheps(t)
    \quad\text{in $L^2(\Omega)$ as $h\to 0$},
\end{equation}
for almost all $t\in [0,T]$.
Let now $h\in \RRR \setminus\{0\}$ be arbitrary. We first observe that
\begin{align}
    \label{EST:DIFF:A}
      &  \iO 
        \partial_j^{+h}\big[A'(\nabla \pheps) \big]
        \cdot \partial_j^{+h} \nabla \pheps
    \notag \\ 
    & \quad = 
    \frac 1{h^2}\iO 
       \big[A'(\nabla \pheps(x+he_j))-A'(\nabla \pheps(x)) \big]
        \cdot \big[\nabla \pheps(x+he_j)-\nabla \pheps(x)\big] \,\mathrm dx
    \notag\\ 
    & \quad = 
    \frac 1{h^2}\iO 
       A'(\nabla \pheps(x+he_j) )
        \cdot \big[\nabla \pheps(x+he_j)-\nabla \pheps(x)\big] \,\mathrm dx
    \notag \\ 
    & \qquad 
       -  \frac 1{h^2}\iO 
       A'(\nabla \pheps(x)) 
        \cdot \big[\nabla \pheps(x+he_j)-\nabla \pheps(x)\big] \,\mathrm dx
    \notag \\ 
    & \quad = 
    - \frac 1{h^2}\iO 
       A'(\nabla \pheps(x) )
        \cdot \big[\nabla \pheps(x+he_j)-2\nabla \pheps(x) + \nabla \pheps(x-he_j)\big] \,\mathrm dx
    \notag \\ 
    &\quad  = - 
    \iO 
    A'(\nabla \pheps) 
    \cdot \partial_j^{-h} \partial_j^{+h} \nabla  \pheps
\end{align}
holds almost everywhere on $[0,T]$. Here, the third equality is obtained by a change of variables. By means of identity \eqref{EST:DIFF:A}, we now deduce
\allowdisplaybreaks
\begin{align}
    \label{EST:DIFF:PHI}
    & c_A \sum_{i,j=1}^d \iO \big| \partial_j^{+h} \partial_i \pheps \big|^2 
    = c_A \sum_{j=1}^d \iO \big| \partial_j^{+h} \nabla \pheps \big|^2 
    \notag\\
    &\quad
    = c_A \sum_{j=1}^d \iO \frac{1}{h^2} \big| \nabla \pheps(x+he_j) - \nabla \pheps (x)\big|^2 \,dx
    \notag\\
    &\quad
    \le \sum_{j=1}^d \iO \frac{1}{h^2} 
        \big[ A'(\nabla\pheps(x+he_j)) - A'(\nabla \pheps(x)) \big] 
        \cdot \big[ \nabla \pheps(x+he_j) - \nabla \pheps(x) \big] \,dx
    \notag\\
    &\quad
    = \sum_{j=1}^d \iO 
        \partial_j^{+h}\big[A'(\nabla \pheps) \big]
        \cdot \partial_j^{+h} \nabla \pheps
    = - \sum_{j=1}^d \iO 
        A'(\nabla \pheps) 
        \cdot \partial_j^{-h} \partial_j^{+h} \nabla  \pheps 
    \notag\\
    &\quad
    = - \sum_{j=1}^d \iO 
        A'(\nabla \pheps) 
        \cdot \nabla (\partial_j^{-h}\partial_j^{+h} \pheps )
    =  \iO 
        \div (A'(\nabla \pheps) )
        \sum_{j=1}^d \partial_j^{-h} \partial_j^{+h} \pheps   
\end{align}
\allowdisplaybreaks[0]%
almost everywhere on $[0,T]$.
Here, from the second to the third line, we used the monotonicity of $A'$ in~\ref{ass:anisotropy}, and the last equality is obtained via integration by parts and the linearity of the integral.
Passing to the limit $h\to 0$, and invoking the convergences \eqref{CONV:PHEPSJI} and \eqref{CONV:PHEPSJJ}, we thus conclude
\begin{align}
    \label{EST:H2}
    c_A \iO | D^2 \pheps|^2 
    \le \iO \div (A'(\nabla \pheps) )
        \Delta \pheps 
\end{align}
almost everywhere on $[0,T]$. This estimate can be used to estimate the second integral on the left-hand side of \eqref{EQ:PHI}. We further recall that due to \eqref{pot:approx}, the third integral on the left-hand side of \eqref{EQ:PHI} is nonnegative. Hence, we infer 
\begin{align}
    \label{EQ:PHI:2}
	& \iO \Phi_\delta(\pheps(t))
	+ c_A \int_0^t\iO | D^2 \pheps|^2
	\le
	\iO \Phi_\delta(\ph_{0,\delta})
	- \int_0^t\iO \psi_2''(\pheps)|\nabla \pheps|^2
\end{align}
for all $t\in [0,T]$.

Using \eqref{UNI:1}, the second integral on the right-hand side of \eqref{EQ:PHI:2} can be bounded by
\begin{equation}
    \label{EST:PHI:1}
    \left| \int_0^t\iO \psi_2''(\pheps)|\nabla \pheps|^2 \, \right|
    \le \norma{\psi_2''}_{L^\infty([-1,1])} \norma{\nabla \phieps}_{\L2 {\Lx2}}^2
    \le C
\end{equation}
for all $t\in [0,T]$.
In view of \eqref{prop:Phi}, we know that $\Phi$ is convex with $\Phi(0)=0$. Since $|\phi_{0,\delta}| \le 1-\delta$ almost everywhere in $\Omega$, we thus obtain
\begin{align}
    \Phi_\delta(\ph_{0,\delta})
    = \Phi(\ph_{0,\delta})
    = \Phi\big((1-\delta)\ph_{0}\big)
    \le (1-\delta) \Phi(\ph_{0})
    \le \Phi(\ph_{0})
\end{align}
almost everywhere in $\Omega$.
Due to condition \eqref{ass:init:data:2}, we thus have
\begin{equation}
    \label{EST:PHI:2}
    \iO \Phi_\delta(\ph_{0,\delta}) \le \iO \Phi(\ph_{0}) \le C.
\end{equation}
Combining \eqref{EQ:PHI:2}, \eqref{EST:PHI:1} and \eqref{EST:PHI:2}, we eventually conclude that
\begin{align}
    \label{UNI:D2PHE}
	\norma{\Phi_\delta(\pheps)}_{\L\infty{\Lx1}}
	   + \norma{D^2\pheps}_{\L2 {\Lx2}} \leq C.
\end{align}
This bound on $\Phi_\delta(\pheps)$ is essential, as it can be used to control deviations of the values of $\pheps$ from the physical range $[-1,1]$ in a quantitative way depending on the size of $\delta$.
This behavior is exactly what we expect, as the mobility degenerates as $\delta$ approaches zero.
A Taylor expansion of $\Phi$ at the point $1-\delta$ yields
\begin{align}
        \Phi_\delta(z) 
        &= \Phi_\delta(1-\delta) 
            + \Phi_\delta'(1-\delta)(z-1+\delta)
            + \frac 12 \Phi_\delta''(\tau)(z-1+\delta)^2      
\end{align}
for all $z> 1$ and some $\tau\in [1-\delta,z]$. Recalling \ref{ass:mob} as well as
\begin{align*}
    &\Phi_\delta(1-\delta) = \Phi(1-\delta) \ge 0,
    \\
    &\Phi_\delta'(1-\delta) = \Phi'(1-\delta) \ge 0,
    \\
    &\Phi_\delta''(\tau) = \frac{1}{\mobe(\tau)} 
        = \frac{1}{\mob(1-\delta)},
\end{align*}
we infer
\begin{align*}
    \Phi_\delta(z) 
    &\ge \frac{1}{2\,\mob(1-\delta)}(z-1+\delta)^2  
    \ge \frac{1}{2\,\big[1-(1-\delta)^2\big]^m B^* }(z-1)^2  
    \notag \\
    &\ge \frac{1}{2\,\big[2\delta\big]^m B^* }(z-1)^2 
    = \frac{1}{2^{m+1}\,\delta^m\, B^*}(|z|-1)^2 
\end{align*}
for all $z> 1$.
Similarly, by performing a Taylor expansion of $\Phi_\delta$ at the point $-1+\delta$, we deduce
\begin{align*}
    \Phi_\delta(z) \ge \frac{1}{2^{m+1}\,\delta^m\, B^*}(|z|-1)^2
\end{align*}
for all $z< -1$. In summary, this shows that
\begin{align*}
    (|z|-1)_+^2 \le 2^{m+1}\,\delta^m\, B^*\,\Phi_\delta(z)
    \quad\text{for all $z\in\RRR$},
\end{align*}
where the subscript ``$+$'' denotes the positive part of the corresponding function.
In combination with \eqref{UNI:D2PHE}, we thus conclude
\begin{align}
    \label{sup:pheps}
    \norma{(|\pheps|-1)_+}_{\L\infty {\Lx2}}^2 =
	\esssup_{t \in [0,T]} \iO \big(|\pheps(t)|-1\big)_+^2 \leq C \delta ^m.
\end{align}


To derive a uniform bound on $w_\delta$, we recall once more that due to \eqref{REG:DNG:E}, we have
\begin{equation}
\label{REG:APPHE}
    \pheps(t)\in H^2(\Omega)
    \quad\text{and}\quad
    A'(\nabla \pheps(t)) \in {\Hx1}
\end{equation}
for almost all $t\in [0,T]$.
Hence, Lemma~\ref{LEM:LIPSCHITZ} implies that 
\begin{align*}
    \big\| D\big[ A'\big(\nabla\pheps \big) \big] \big\|_{\Lx2} 
    \le 
    C \norma{ D^2 \pheps}_{\Lx2}  
\end{align*}
almost everywhere on $[0,T]$.
Eventually, using the uniform bound \eqref{UNI:D2PHE}, we conclude
\begin{align}
    \label{UNI:WE}
    &\norma{w_\delta}_{\L2 {\Lx2}} 
    = \norma{\div \big[A'(\nabla\pheps)\big]}_{\L2 {\Lx2}}
    \notag \\
    &\quad 
    \le C \norma{D \big[A'(\nabla\pheps)\big]}_{\L2 {\Lx2}}
    \le C \norma{D^2\pheps}_{\L2 {\Lx2}} \le C.
\end{align}

In summary, collecting the above estimates \eqref{EST:MOBE}, \eqref{EST:DMOBE}, \eqref{UNI:PHE}, 
\eqref{UNI:JE}, \eqref{UNI:DTPHE}, \eqref{UNI:D2PHE} and \eqref{UNI:WE}, we have shown 
\begin{align}
    \label{UNI:ALL:1}
    & 
    \norma{\pheps}_{\H1 {{\Hx1}^*}\cap \L\infty {\Hx1} \cap \L2 {\Hx2}}
    + \norma{J_\delta}_{\L2 {{L^2(\Omega;\RRR^d)}}}
     \notag \\ 
    & \quad 
    + \norma{\mobe(\pheps)}_{L^\infty(Q) \cap \L2{\Hx1}}
    + \norma{w_\delta}_{\L2 {\Lx2}}
    \le C,
\end{align}
as well as
\begin{align}
    \label{UNI:ALL:2}
    \norma{(|\pheps|-1)_+}_{\L\infty {\Lx2}} \le C \delta^{m/2}.
\end{align}

\color{black}


\subsection{Passing to the limit}
We now exploit the uniform estimates \eqref{UNI:ALL:1} and \eqref{UNI:ALL:2} to pass to the limit $\delta \to 0$.
Based on \eqref{UNI:ALL:1}, the Banach--Alaoglu theorem implies that
there exist suitable functions $\ph$, $J$, $w$ and $b_*$ such that
\begin{alignat}{2}
    \label{CONV:PH}
    &\pheps \to \ph \quad 
    &&\text{weakly-star in $\L\infty {\Hx1}$},
    \notag \\
    &{}&&\text{and weakly in $\H1 \Vp \cap \L2 {\Hx2}$},
    \\
    \label{CONV:J}
    &J_\delta \to J \quad  
    &&\text{weakly in $\L2 {L^2(\Omega;\RRR^d)}$},
     \\
    &
    \label{CONV:B*}
    \mobe(\pheps) \to b_* \quad 
    &&\text{weakly-star in $L^\infty(Q)$ and weakly in $\L2 {\Hx1}$},
    \\
    &
    \label{CONV:W}
    w_\delta \to w \quad  
    &&\text{weakly in $\L2 {L^2(\Omega)}$},
\end{alignat}
as $\delta \to 0$ along a non-relabeled subsequence.
Moreover, the Aubin--Lions lemma (see, e.g., \cite[Sect.~8, Cor.~4]{Simon}) then yields
\begin{align}
    \label{CONV:PH:S}
    \pheps \to \ph \quad \text{strongly in $\C0 {\Lx2}$, and $\aeQ$,}
\end{align}
up to subsequence extraction. In particular, this entails that
\begin{align*}
    \ph\vert_{t=0} = \underset{\delta\to 0}{\lim}\,\pheps\vert_{t=0} = \phi_0
    \quad\text{\aeO.}
\end{align*}
By passing to the limit $\delta \to 0$ in \eqref{UNI:ALL:2}, we further infer that
\begin{align*}
     |\ph| \leq 1 \quad \aeQ.
\end{align*}

To pass to the limit in the anisotropic term, the weak convergence
$ \nabla \pheps \to \nabla \ph$ in $\L2 {L^2(\Omega;\RRR^d)}$ as $\delta \to 0$ does not suffice. 
However, due to the strong monotonicity of $A'$ (cf.~\ref{ass:anisotropy}), we even have
\begin{align*}
   &c_A \norma{\nabla \pheps-\nabla \ph}_{L^2(\Omega;\RRR^d)}^2
   \leq \iO \big(A'(\nabla \pheps)-A'(\nabla \ph)\big) \cdot \nabla (\pheps- \ph)
   \\ 
   &\quad 
   = \iO w_\delta (\pheps- \ph)
   - \iO A'(\nabla \ph) \cdot\nabla (\pheps- \ph),
\end{align*}
and the terms on the right-hand side converge to zero as $\delta \to 0$. This directly entails
\begin{align}
    \label{CONV:GPH:S}
    \nabla \pheps \to \nabla \ph \quad \text{strongly in $\L2 {L^2(\Omega;\RRR^d)}$ as $\delta \to 0$.} 
\end{align}

Due to the growth condition on $A'$ in \ref{ass:anisotropy}, we infer 
\begin{align}
    \label{CONV:A'}
    A' (\nabla \pheps) &\to A' (\nabla \ph)
    \quad \text{strongly in $L^2(0,T; L^2(\Omega;\RRR^d))$},
\end{align}
as $\delta\to 0$, by Lebesgue's general convergence theorem (see, e.g., \cite[Section~3.25]{Alt}).

Next, we want to identify $b_*$ with $b(\phi)$. 
Recalling that $|\ph|\le 1$ \aeQ, that $\mobe \to b$ as $\delta\to 0$ pointwise on $\RRR$, and that $\mobe$ is Lipschitz continuous with a Lipschitz constant independent of $\delta$, we deduce that, for almost all $(x,t) \in Q$,
\begin{align*}
    &\big|\mobe\big(\pheps(x,t)\big) - \mob\big(\ph(x,t)\big)\big|
    \\[0.5ex]
    &\quad
    \le |\mobe\big(\pheps(x,t)\big) - \mobe\big(\ph(x,t)\big)|
        + |\mobe\big(\ph(x,t)\big) - \mob\big(\ph(x,t)\big)|
    \\[0.5ex]
    &\quad
    \le C|\pheps(x,t) - \ph(x,t)| 
        + |\mobe\big(\ph(x,t)\big) - \mob\big(\ph(x,t)\big)|
    \to 0
\end{align*}
as $\delta\to 0$.
As the pointwise limit and weak limit in \eqref{CONV:B*} coincide (see, e.g., \cite[Corollary~9.2c]{DiBenedetto}), we thus conclude that
\begin{align*}
    b_* = b(\ph) \quad\aeQ.
\end{align*}
Consequently, as $\delta\to 0$, we have
\begin{align}
    \label{CONV:B}
    \mobe(\pheps) \to \mob(\ph) \quad  
    &\text{weakly-star in $L^\infty(Q)$, weakly in $\L2 {\Hx1}$},
    \notag\\
    &\text{and \aeQ.}
\end{align}
Collecting the above information, we have thus shown that the triplet $(\varphi,J,w)$ fulfills all the conditions demanded in item \ref{THM:MAIN:1} of Theorem~\ref{THM:MAIN}.

Our final goal is to pass to the limit in the approximate weak formulation \eqref{WF:NDG:ALT}. 
This can essentially be done by following the line of argument in \cite{EG}. 
To this end, let $\bet\in H^1(\Omega;\RRR^d) \cap L^\infty(\Omega;\RRR^d)$ be arbitrary. We first recall that 
\begin{align}
\label{CR:DIVB}
    \div\big( \mobe(\pheps) \bet \big) 
    = \mobe'(\pheps) \nabla\pheps \cdot \eta
        + \mobe(\pheps) \div(\bet)
    \quad\aeQ.
\end{align}
Since $\mobe(\pheps)$ is bounded uniformly in $\delta$ (cf.~\eqref{EST:MOBE}) and $\div(\bet)\in L^2(\Omega)$, \eqref{CONV:B} along with Lebesgue's dominated convergence theorem implies, as $\delta\to 0$,
\begin{align}
    \label{CONV:B:L2}
    \mobe(\pheps) \div(\bet) \to \mob(\ph) \div(\bet)
    \quad\text{strongly in $L^2(0,T; L^2(\Omega;\RRR^d))$.}
\end{align}

We next claim that, as $\delta\to 0$,
\begin{align}
    \label{CONV:DIVB:L2}
    \mobe'(\pheps) \nabla\pheps 
    \to \mob'(\ph) \nabla\ph 
    \quad\text{strongly in $L^2(0,T; L^2(\Omega;\RRR^d))$.}
\end{align}

To show this, we recall that $|\ph|\le 1$ \aeQ, and employ the decomposition
\begin{align}
    \label{DECOMP}
    &\int_Q \big| \mobe'(\pheps) \nabla\pheps - \mob'(\ph) \nabla\ph\big|^2
    \notag \\
    &\quad =\int\limits_{\{|\ph|=1\}} \big| \mobe'(\pheps) \nabla\pheps - \mob'(\ph) \nabla\ph\big|^2
        + \int\limits_{\{|\ph|<1\}} \big| \mobe'(\pheps) \nabla\pheps - \mob'(\ph) \nabla\ph\big|^2.
\end{align}
Here, we use the notation $\{|\ph|=1\}: = \{ (x,t) \in Q : |\ph(x,t) |=1\}$.
Since $\nabla\ph = 0$ almost everywhere on $\{|\ph|=1\}$ (see, e.g., \cite[§7.4,~Lemma~7.7]{Gilbarg2001}), we use the strong convergence \eqref{CONV:GPH:S} along with the fact that $\mobe'$ is bounded uniformly with respect to $\delta$ (see \eqref{EST:DMOBE}) to deduce
\begin{align}
    \label{CONV:T1}
    &\int\limits_{\{|\ph|=1\}} \big| \mobe'(\pheps) \nabla\pheps - \mob'(\ph) \nabla\ph\big|^2
    = \int\limits_{\{|\ph|=1\}} \big| \mobe'(\pheps) \nabla\pheps\big|^2
    \notag \\
    &\quad
    \le C \int\limits_{\{|\ph|=1\}} \big|\nabla\pheps\big|^2
    \to C \int\limits_{\{|\ph|=1\}} \big|\nabla\ph\big|^2 = 0
\end{align}
as $\delta\to 0$.
To estimate the second integral on the right-hand side of \eqref{DECOMP}, we recall from \eqref{CONV:PH:S} that there exists a null set $\mathcal N\subset Q$ such that
\begin{align}
    \label{CONV:PHE:PTW}
    \pheps(x,t) \to \ph(x,t) 
    \quad\text{for all $(x,t) \in Q\setminus\mathcal N$.}
\end{align}
Let us now fix an arbitrary $(x,t) \in \{|\ph|<1\} \setminus\mathcal N$. Then, there exists $\alpha>0$ such that 
\begin{equation*}
    |\ph(x,t)| \le 1 - 2\alpha .
\end{equation*}
If $\delta$ is chosen sufficiently small, it further holds that
\begin{equation*}
    |\pheps(x,t) - \ph(x,t)| \le \alpha
\end{equation*}
because of \eqref{CONV:PHE:PTW}. Choosing $\delta$ possibly even smaller such that $\delta<\alpha$, we have
\begin{equation*}
    |\pheps(x,t)| \le |\pheps(x,t) - \ph(x,t)| + |\ph(x,t)| < 1-\alpha < 1-\delta,
\end{equation*}
and thus
\begin{equation*}
    \mobe'\big(\pheps(x,t)\big) 
    = \mob'\big(\pheps(x,t)\big) 
    \to \mob'\big(\ph(x,t)\big) 
    \quad\text{as $\delta\to 0$}
\end{equation*}
due to the continuity of $\mob'$ on $(-1,1)$. This proves that
\begin{equation*}
    \mobe'(\pheps) \to \mob'(\ph)
    \quad\text{a.e.~in $\{|\ph|<1\}$ as $\delta\to 0$.}
\end{equation*}
Recalling \eqref{CONV:GPH:S} and that $\mobe'$ is bounded uniformly in $\delta$ (see \eqref{EST:DMOBE}), we thus conclude
\begin{align}
    \label{CONV:T2}
    \int\limits_{\{|\ph|<1\}} \big| \mobe'(\pheps) \nabla\pheps - \mob'(\ph) \nabla\ph\big|^2
    \to 0
\end{align}
as $\delta\to 0$. Using \eqref{CONV:T1} and \eqref{CONV:T2} to pass to the limit in the right-hand side of \eqref{DECOMP}, we have verified the convergence claimed in \eqref{CONV:DIVB:L2}.

Eventually, the strong convergences \eqref{CONV:B:L2} and \eqref{CONV:DIVB:L2} along with the decomposition \eqref{CR:DIVB} show that, as $\delta\to 0$,
\begin{align}
    \label{CONV:BET}
    \div\big( \mobe(\pheps) \bet \big) 
    \to \div\big( \mob(\ph) \bet \big) 
\end{align}
for every $\bet\in H^1(\Omega;\RRR^d) \cap L^\infty(\Omega;\RRR^d)$.

In order to pass to the limit in \eqref{WF:NDG:ALT}, we still need to show that, as $\delta\to 0$,
\begin{align}
    \label{CONV:BETPSI}
    \mobe(\pheps) \psi_\delta''(\pheps) \nabla\pheps
    \to \mob(\ph) \psi''(\ph) \nabla\phi
    \quad\text{strongly in $L^2(0,T; L^2(\Omega;\RRR^d))$.}
\end{align}

This convergence can be shown similarly as the one in \eqref{CONV:DIVB:L2}. 
Recalling \eqref{ID:BPSI}, \eqref{mob:approx} and \eqref{DEF:POT:APPROX}, we first notice that for all $r\in [-1,1]$,
\begin{align}
    \label{EST:MPSI}
    &| \mobe(r) \psi_\delta''(r) |
    \le \| b\psi'' \|_{L^\infty([-1,1])}
    = \| FB + b\psi_2'' \|_{L^\infty([-1,1])}
    \le C.
\end{align}
Moreover, we have the pointwise convergence
\begin{equation*}
    \mobe(r) \psi_\delta''(r) 
    \to \mob(r) \psi''(r)
    \quad\text{as $\delta\to 0$}
\end{equation*}
for all $r\in [-1,1]$.
Since $|\ph|\le 1$ \aeQ, we make the decomposition
\begin{align}
    \label{DECOMP:2}
    &\int_Q \big| \mobe(\pheps) \psi_\delta''(\pheps) \nabla\pheps
        - \mob(\ph) \psi''(\ph) \nabla\phi\big|^2
    \notag\\
    &\quad= \int\limits_{\{|\ph|=1\}} \big| \mobe(\pheps) \psi_\delta''(\pheps) \nabla\pheps
        - \mob(\ph) \psi''(\ph) \nabla\phi\big|^2
    \notag\\
    &\qquad    
    + \int\limits_{\{|\ph|<1\}} \big| \mobe(\pheps) \psi_\delta''(\pheps) \nabla\pheps
        - \mob(\ph) \psi''(\ph) \nabla\phi\big|^2.
\end{align}
Since $\nabla\ph = 0$ almost everywhere on $\{|\ph|=1\}$ (see, e.g., \cite[§7.4,~Lemma~7.7]{Gilbarg2001}), we use the strong convergence \eqref{CONV:GPH:S} along with the uniform bound \eqref{EST:MPSI} to deduce
\begin{align}
    \label{CONV:I1}
    &\int\limits_{\{|\ph|=1\}} \big| \mobe(\pheps) \psi_\delta''(\pheps) \nabla\pheps
        - \mob(\ph) \psi''(\ph) \nabla\phi\big|^2
    \notag\\
    &\quad
    = \int\limits_{\{|\ph|=1\}} \big| \mobe(\pheps) \psi_\delta''(\pheps) \nabla\pheps \big|^2
    \le C \int\limits_{\{|\ph|=1\}} \big| \nabla\pheps \big|^2
    \to C \int\limits_{\{|\ph|=1\}} \big| \nabla\ph \big|^2
    = 0
\end{align}
as $\delta\to 0$.
Proceeding as in the proof of \eqref{CONV:DIVB:L2}, we further use the pointwise convergence \eqref{CONV:PHE:PTW} to show that for all $(x,t) \in \{|\ph|<1\} \setminus\mathcal N$, it holds
\begin{equation*}
    |\pheps(x,t)|\le 1-\delta
    \quad\text{and thus,}\quad
    \mobe(\pheps) \psi_\delta''(\pheps) = \mob(\pheps) \psi''(\pheps)
\end{equation*}
if $\delta$ is chosen sufficiently small.
Invoking the continuity of $\mob \psi''$ on $(-1,1)$, we then deduce that
\begin{equation*}
    \mobe(\pheps) \psi_\delta''(\pheps) 
    \to \mob(\ph) \psi''(\ph)
    \quad\text{a.e.~in $\{|\ph|<1\}$ as $\delta\to 0$.}
\end{equation*}
Hence, employing Lebesgue's dominated convergence theorem along with the uniform bound \eqref{EST:MPSI}, we infer
\begin{equation}
    \label{CONV:I2}
    \int\limits_{\{|\ph|<1\}} \big| \mobe(\pheps) \psi_\delta''(\pheps) \nabla\pheps
        - \mob(\ph) \psi''(\ph) \nabla\phi\big|^2
    \to 0 \quad\text{as $\delta\to 0$}.
\end{equation}
Eventually, combining \eqref{CONV:I1} and \eqref{CONV:I2}, we pass to the limit in the right-hand side of \eqref{DECOMP:2}. This verifies \eqref{CONV:BETPSI}.

Finally, by combining the convergences \eqref{CONV:PH}--\eqref{CONV:A'}, \eqref{CONV:BET} and \eqref{CONV:BETPSI}, we can pass to the limit $\delta\to 0$ in the approximate weak formulation \eqref{WF:NDG:ALT}. This means that item \ref{THM:MAIN:2} of Theorem~\ref{THM:MAIN} is verified and thus, the proof is complete.

\section*{Appendix: Compositions of Lipschitz and Sobolev functions}
\renewcommand\thesection{A}
\setcounter{theorem}{0}
\setcounter{equation}{0}

As a supplement, we present the following elementary, but very useful lemma on compositions of Lipschitz and Sobolev functions.

\begin{lemma} \label{LEM:LIPSCHITZ}
    Let $\Omega$ be a $d$-dimensional flat torus and let $1\le p<\infty$.
    Let $n,m \in \mathbb{N}$, let ${v:\RRR^n \to \RRR^m}$ be Lipschitz continuous with Lipschitz constant $L\ge 0$, and let ${u\in W^{1,p}(\Omega;\RRR^n)}$. We further assume that ${v\circ u} \in L^p(\Omega;\RRR^m)$. Then, ${v\circ u}:\Omega\to\RRR^m$ is weakly differentiable and there exists a constant $C\ge 0$ depending only on $d$, $n$ and $m$ such that the weak derivative satisfies
    \begin{align*}
        \norma{D(v\circ u)}_{L^p(\Omega;\RRR^{m\times d})} \le CL\, \norma{Du}_{L^p(\Omega;\RRR^{n\times d})}.
    \end{align*}
\end{lemma}

\begin{proof}
    Let $i\in \{1,...,m\}$ be arbitrary and let $C\ge 0$ denote a generic constant depending only on $d$, $n$ and $m$, which may change its value from line to line. Then, for almost all $x\in \Omega$, it holds
    \begin{align*}
        &\left|\frac{(v\circ u)_i(x+he_j) - (v\circ u)_i(x)}{h}\right|^2
        \le \frac{\big| (v\circ u)(x+he_j) - (v\circ u)(x) \big|^2}{|h|^2}
        \\
        &\quad
        \le L^2 \; \frac{\big| u(x+he_j) - u(x) \big|^2}{|h|^2}
        = L^2 \; \sum_{k=1}^n \left| \frac{u_k(x+he_j) - u_k(x) }{h} \right|^2
        \\
        &\quad
        \le C L^2\; \left( \sum_{k=1}^n \left| \frac{u_k(x+he_j) - u_k(x) }{h} \right| \right)^2,
    \end{align*}
    where $e_j$ denotes the $j$-th unit vector in $\RRR^d$ and $h\in \RRR\setminus\{ 0\}$.
    This directly entails
    \begin{align*}
        &\left\|\frac{(v\circ u)_i(x+he_j) - (v\circ u)_i(x)}{h}\right\|_{L^p(\Omega)}
        \le C L\; \sum_{k=1}^n \left\| \frac{u_k(x+he_j) - u_k(x) }{h} \right\|_{L^p(\Omega)}
        \\
        &\quad
        \le\, CL\; \sum_{k=1}^n \norma{ \partial_j u_k}_{L^p(\Omega)}
        \le CL\; \norma{Du}_{L^p(\Omega;\RRR^{n\times d})}.
    \end{align*}
    Therefore, for every $i\in\{1,...,n\}$, the difference quotients of $(v\circ u)_i$ are bounded in $L^p(\Omega)$. We conclude (see, e.g.,
    \cite[§5.8.2,~Theorem~3]{Evans2010}) that $v\circ u$ is weakly differentiable with
    \begin{equation*}
        \norma{D(v\circ u)}_{L^p(\Omega;\RRR^{n\times d})} \le C L \norma{Du}_{L^p(\Omega;\RRR^{n\times d})}.
    \end{equation*}
    Since, $v\circ u \in L^p(\Omega;\RRR^n)$, it directly follows that $v\circ u \in W^{1,p}(\Omega;\RRR^n)$ and thus, the proof is complete.
\end{proof}

\section*{Acknowledgements}
\noindent
HG and PK were partially supported by the ``RTG 2339'', which is funded by the Deutsche Forschungsgemeinschaft (DFG, German Research Foundation). The support is gratefully acknowledged.
AS gratefully acknowledges some support 
from the MIUR-PRIN Grant 2020F3NCPX ``Mathematics for industry 4.0 (Math4I4)'', from ``MUR GRANT Dipartimento di Eccellenza'' 2023-2027 and from the Alexander von Humboldt Foundation.
Additionally, AS appreciates his 
affiliation with GNAMPA (Gruppo Nazionale per l'Analisi Matematica, la Probabilità e le loro Applicazioni) of INdAM (Istituto Nazionale di Alta Matematica).

\footnotesize

\bibliography{GKS}

\End{document}
